# Variations on the Bloch-Ogus Theorem

I. Panin, K. Zainoulline

11/03/2002

**Abstract**

In the present paper we discuss questions concerning the arithmetic resolution for étale cohomology. In particular, consider a smooth quasi-projective variety $X$ over a field $k$ together with the local scheme $\mathcal{U} = \operatorname{Spec} \mathcal{O}_{X,x}$ at a point $x \in X$. Let $p: Y \to \mathcal{U}$ be a smooth projective morphism and let $Y_{k(u)}$ denote its fiber over the generic point of a subvariety $u$ of $\mathcal{U}$. We prove there is a Gersten-type exact sequence

$$0 \to H^q_{\text{ét}}(Y, F) \to H^q_{\text{ét}}(Y_{k(\mathcal{U})}, F) \to \coprod_{u \in \mathcal{U}^{(1)}} H^{q-1}_{\text{ét}}(Y_{k(u)}, F(-1)) \to$$

of étale cohomology with coefficients in a locally constant étale sheaf $F$ of $\mathbb{Z}/n\mathbb{Z}$-modules on $Y$ which has finite stalks and $(n, char(k)) = 1$.

## 1 Introduction

In the present paper we discuss questions concerning the arithmetic resolution for étale cohomology.

Let $X$ be a smooth variety over a field $k$ and let $x = \{x_1, \ldots, x_m\} \subset X$ be a finite set of points. We denote by $\mathcal{U} = \operatorname{Spec} \mathcal{O}_{X,x}$ the semi-local scheme at $x$. Then $\mathcal{U}$ is a semi-local regular scheme of geometric type over $k$. Consider the sheaf $\mu_n$ of $n$-th roots of unity on the small étale site $X_{et}$, with $n$ prime to the characteristic of the base field $k$.

In the famous paper of S. Bloch and A. Ogus ([1], Example 2.1 and Theorem 4.2) it was proven for all $i \in \mathbb{Z}$ and $q \geq 0$ there is an exact sequence

$$0 \to H^q(\mathcal{U}, \mu_n^{\otimes i}) \to \coprod_{u \in \mathcal{U}^{(0)}} H^q_u(\mathcal{U}, \mu_n^{\otimes i}) \to \coprod_{u \in \mathcal{U}^{(1)}} H^{q+1}_u(\mathcal{U}, \mu_n^{\otimes i}) \to \cdots \quad (\dagger)$$

consisting of étale cohomology with supports, where $\mathcal{U}^{(p)}$ denotes the set of all points of codimension $p$ in $\mathcal{U}$.

The next step was done by O. Gabber in [6]. He proved that the sequence (†) is exact for cohomology with coefficients in any torsion sheaf $F$ on $X_{et}$ that comes from the base field $k$, i.e., $F = p^*G$ for some sheaf $G$ on $(\operatorname{Spec} k)_{et}$ and the structural morphism $p: X \to \operatorname{Spec} k$. The crucial point of the proof was to show that the sheaf $F$ is effaceable ([3], Definition 2.1.1). Then the exactness of (†) follows immediately by trivial reasons ([3], Proposition 2.1.2).

It turned out that the proof of Gabber can be applied to any cohomology theory with supports that satisfies the same formalism as étale cohomology do. This idea was realized in the paper of J.-L. Colliot-Thélène, R. Hoobler and B. Kahn [3]. Namely, they proved that a cohomology theory with support $h^*$ which satisfies some set of axioms (described in [3], section 5.1) is effaceable and, thus, there is an exact sequence similar to (†). In particular, one gets (†) for the case when $\mathcal{U}$ is replaced by the product $\mathcal{U} \times_k T$, where $T$ is a smooth variety over $k$ ([3], Theorem 8.1.1).

In the same paper by using different arguments ([3], Remark 8.1.2.(3), Theorem B.2.1) the exactness of (†) for the case when $\dim \mathcal{U} = 1$ and $\mu_n^{\otimes i}$ is replaced by any bounded complex of locally constant constructible torsion sheaves on $X_{et}$ with torsion prime to the characteristic of the base field $k$ was proven. The goal of this paper is to prove the latter case for any dimension of the scheme $\mathcal{U}$. Namely, we want to prove the following

**1.1. Theorem.** *Let $X$ be a smooth quasi-projective variety over an infinite field $k$. Let $x = \{x_1, \ldots, x_m\} \subset X$ be a finite set of points and $\mathcal{U} = \operatorname{Spec} \mathcal{O}_{X,x}$ be the semi-local scheme at $x$. Let $\mathcal{C}$ be a bounded complex of locally constant constructible sheaves of $\mathbb{Z}/n\mathbb{Z}$-modules on $X_{et}$ with $n$ prime to the characteristic of the base field $k$. Then the $E^1$-terms of the coniveau spectral sequence yield an exact sequence*

$$0 \to H^q(\mathcal{U}, \mathcal{C}) \to \coprod_{u \in \mathcal{U}^{(0)}} H^q_u(\mathcal{U}, \mathcal{C}) \to \coprod_{u \in \mathcal{U}^{(1)}} H^{q+1}_u(\mathcal{U}, \mathcal{C}) \to \cdots$$

*of étale hypercohomology with supports.*

**1.2. Corollary.** *Let $R$ be a semi-local regular ring of geometric type over an infinite field $k$. We denote by $\mathcal{U} = \operatorname{Spec} R$ the respective semi-local affine scheme. Let $\mathcal{C}$ be a bounded complex of locally constant constructible sheaves* of $\mathbb{Z}/n\mathbb{Z}$-modules on $\mathcal{U}_{et}$ with $n$ prime to the characteristic of $k$. Then there is an exact sequence

$$0 \to H^q_{\text{ét}}(\mathcal{U},\mathcal{C}) \to H^q_{\text{ét}}(\operatorname{Spec} k(\mathcal{U}),\mathcal{C}) \to \coprod_{u \in \mathcal{U}^{(1)}} H^{q-1}_{\text{ét}}(\operatorname{Spec} k(u),\mathcal{C}(-1)) \to \cdots$$

where $k(\mathcal{U})$ is the function field of $\mathcal{U}$, $k(u)$ is the residue field of $u$ and $\mathcal{C}(i) = \mathcal{C} \otimes \mu_n^{\otimes i}$.

*Proof.* Follows by purity for étale cohomology (see the arguments on p.36–37 of the proof of (1.4), [3]) applied to the complex from Theorem 1.1. □

**1.3. Corollary.** *Let $R$ be a semi-local regular ring of geometric type over an infinite field $k$. Let $p : Y \to \mathcal{U}$ be a smooth proper morphism, where $\mathcal{U} = \operatorname{Spec} R$. Let $F$ be a locally constant constructible sheaf of $\mathbb{Z}/n\mathbb{Z}$-modules on $Y_{et}$ with $n$ prime to the characteristic of the base field $k$. Then there is an exact sequence*

$$0 \to H^q_{\text{ét}}(Y, F) \to H^q_{\text{ét}}(Y_{k(\mathcal{U})}, F) \to \coprod_{u \in \mathcal{U}^{(1)}} H^{q-1}_{\text{ét}}(Y_{k(u)}, F(-1)) \to \cdots$$

*where $Y_{k(u)} = \operatorname{Spec} k(u) \times_{\mathcal{U}} Y$.*

*Proof.* The cohomology of $Y$ with coefficients in $F$ coincide with the hypercohomology of $\mathcal{U}$ with coefficients in the total direct image $Rp_*F$ (see [7], VI.4.2). By ([10], Theorem A) there exists a bounded complex $\mathcal{C}$ of locally constant constructible sheaves that is quasi-isomorphic to the complex $Rp_*F$. Now replace $Rp_*F$ by $\mathcal{C}$ and apply the previous corollary. We get the desired exact sequence. □

**1.4. Remark.** Corollary 1.3 still holds if one replaces the sheaf $F$ by any bounded complex of locally constant constructible sheaves of $\mathbb{Z}/n\mathbb{Z}$-modules on $Y_{et}$ with $n$ prime to the characteristic of the base field $k$.

**1.5. Remark.** By arguments with transfers one can extend all the results to the case of a finite field $k$.

As in [3] in order to prove Theorem 1.1 we show that any cohomology functor with supports $F$ which satisfies some set of axioms (given in section 2) is effaceable in the sense of Definition 4.1. First, we do this for the case when $F$ is defined over the base field (Theorem 4.2). Then we do this for the general situation, when $F$ is defined over some smooth variety (Theorem 5.4). Finally, we check that the étale cohomology functor $F(X, Z) = H^*_Z(X, \mathcal{C})$ from 1.1 satisfies these axioms and, thus, is effaceable. So that Theorem 1.1 follows from Lemma 5.7 which is a slightly modified version of Proposition 2.1.2, [3].

The key difference from [3] is that we use Geometric Presentation Lemma of Ojanguren and Panin ([8], 10.1) instead of Gabber's. It allows us to use transfer arguments which were motivated by paper [13] and developed in [8], [9] and [14]. The main advantage of this techniques is that it can be extended to the case we are interesting in – when the sheaf of coefficients $F$ is defined not only over the base field $k$ but over some smooth variety over $k$. In order to realise this we also use the main result of paper [10].

## 2   Definitions and Notations

**2.1. Notation.** In the present paper all schemes are assumed to be Noetherian and separated. By $k$ we denote a fixed ground field. A variety over $k$ is an integral scheme of finite type over $k$. To simplify the notation sometimes we will write $k$ instead of the scheme $\operatorname{Spec} k$. We will write $X_1 \times X_2$ for the fibred product $X_1 \times_k X_2$ of two $k$-schemes. By $\mathcal{U}$ we denote a regular semi-local scheme of geometric type over $k$, i.e., $\mathcal{U} = \operatorname{Spec} \mathcal{O}_{X,x}$ for a smooth affine variety $X$ over $k$ and a finite set of points $x = \{x_1, \ldots, x_n\}$ of $X$. By $\mathcal{X}$ we denote a relative curve over $\mathcal{U}$ (see 3.1.(i)). By $\mathcal{Z}$ and $\mathcal{Y}$ we denote closed subsets of $\mathcal{X}$. By $Z$ and $Z'$ we denote closed subsets of $\mathcal{U}$. Observe that $\mathcal{X}$ and $\mathcal{U}$ are essentially smooth over $k$ and all schemes $\mathcal{X}, \mathcal{Z}, \mathcal{Y}, Z, Z'$ are of finite type over $\mathcal{U}$.

**2.2. Notation.** Let $U$ be a $k$-scheme. Denote by $Cp(U)$ a category whose objects are couples $(X, Z)$ consisting of an $U$-scheme $X$ of finite type over $U$ and a closed subset $Z$ of the scheme $X$ (we assume the empty set is a closed subset of $X$). Morphisms from $(X, Z)$ to $(X', Z')$ are those morphisms $f : X \to X'$ of $U$-schemes that satisfy the property $f^{-1}(Z') \subset Z$. The composite of $f$ and $g$ is $g \circ f$.

**2.3. Notation.** Denote by $F : Cp(U) \to Ab$ a contravariant functor from the category of couples $Cp(U)$ to the category of (graded) abelian groups. We shall write $F_Z(X)$ for $F(X, Z)$ having in mind the notation used for cohomology with supports.

Now notions of a homotopy invariant functor, a functor with transfers and a functor that satisfies vanishing property will be given.

**2.4. Definition.** A contravariant functor $F : Cp(U) \to Ab$ is said to be homotopy invariant if for each $U$-scheme $X$ smooth or essentially smooth over $k$ and for each closed subset $Z$ of $X$ the map $F_Z(X) \to F_{Z \times \mathbb{A}^1}(X \times \mathbb{A}^1)$ induced by the projection $X \times \mathbb{A}^1 \to X$ is an isomorphism.

**2.5. Definition.** One says a contravariant functor $F : Cp(U) \to Ab$ satisfies vanishing property if for each $U$-scheme $X$ one has $F(X, \emptyset) = 0$.

**2.6. Definition.** A contravariant functor $F : Cp(U) \to Ab$ is said to be endowed with transfers if for each finite flat morphism $\pi : X' \to X$ of $U$-schemes and for each closed subset $Z \subset X$ it is given a homomorphism of abelian groups $\operatorname{Tr}_X^{X'} : F_{\pi^{-1}(Z)}(X') \to F_Z(X)$ and the family $\{\operatorname{Tr}_X^{X'}\}$ satisfies the following properties:

(i) for each fibred product diagram of $U$-schemes with a finite flat morphism $\pi$

$$\begin{array}{ccc} X' & \xleftarrow{f'} & X'_1 \\ \pi \downarrow & & \downarrow \pi_1 \\ X & \xleftarrow{f} & X_1 \end{array}$$

and for each closed subset $Z \subset X$ the diagram

$$\begin{array}{ccc} F_{Z'}(X') & \xrightarrow{F(f_1)} & F_{Z'_1}(X'_1) \\ \operatorname{Tr}_X^{X'} \downarrow & & \downarrow \operatorname{Tr}_{X_1}^{X'_1} \\ F_Z(X) & \xrightarrow{F(f)} & F_{Z_1}(X_1) \end{array}$$

is commutative, where $Z' = \pi^{-1}(Z)$, $Z_1 = f^{-1}(Z)$ and $Z'_1 = \pi_1^{-1}(Z_1)$;

(ii) if $\pi : X'_1 \amalg X'_2 \to X$ is a finite flat morphism of $U$-schemes, then for each closed subset $Z \subset X$ the diagram

$$\begin{array}{ccc} F_{Z'}(X'_1 \amalg X'_2) & \xrightarrow{\text{``+''}} & F_{Z'_1}(X'_1) \oplus F_{Z'_2}(X'_2) \\ & \searrow{\scriptstyle \text{Tr}_X^{X'_1 \amalg X'_2}} & \downarrow{\scriptstyle \text{Tr}_{X_1}^{X'_1} + \text{Tr}_{X_2}^{X'_2}} \\ & & F_Z(X) \end{array}$$

is commutative, where $Z' = \pi^{-1}(Z)$, $Z'_1 = X'_1 \cap Z'$ and $Z'_2 = X'_2 \cap Z'$;

(iii) if $\pi : (X', Z') \to (X, Z)$ is an isomorphism in $Cp(U)$, then two maps $\text{Tr}_X^{X'}$ and $F(\pi)$ are inverses of each other, i.e.

$$F(\pi) \circ \text{Tr}_X^{X'} = \text{Tr}_X^{X'} \circ F(\pi) = \text{id}.$$

# 3 The Specialization Lemma

The following definition is inspired by the notion of a good triple used by Voevodsky in [13].

**3.1. Definition.** Let $\mathcal{U}$ be a regular semi-local scheme of geometric type over the field $k$. A triple $(\mathcal{X}, \delta, \mathfrak{f})$ consisting of an $\mathcal{U}$-scheme $p : \mathcal{X} \to \mathcal{U}$, a section $\delta : \mathcal{U} \to \mathcal{X}$ of the morphism $p$ and a regular function $\mathfrak{f} \in \Gamma(\mathcal{X}, \mathcal{O}_\mathcal{X})$ is called a perfect triple over $\mathcal{U}$ if $\mathcal{X}$, $\delta$ and $\mathfrak{f}$ satisfy the following conditions:

(i) the morphism $p$ can be factorized as $p : \mathcal{X} \xrightarrow{\pi} \mathbb{A}^1 \times \mathcal{U} \xrightarrow{pr} \mathcal{U}$, where $\pi$ is a finite surjective morphism and $pr$ is the canonical projection on the second factor;

(ii) the vanishing locus of the function $\mathfrak{f}$ is finite over $\mathcal{U}$;

(iii) the scheme $\mathcal{X}$ is essentially smooth over $k$ and the morphism $p$ is smooth along $\delta(\mathcal{U})$;

(iv) the scheme $\mathcal{X}$ is irreducible.

**3.2. Remark.** The property (i) says that $\mathcal{X}$ is an affine curve over $\mathcal{U}$. The property (iii) implies that $\mathcal{X}$ is a regular scheme. Since $\mathcal{X}$ and $\mathbb{A}^1 \times \mathcal{U}$ are regular schemes by ([5], 18.17) the morphism $\pi : \mathcal{X} \to \mathbb{A}^1 \times \mathcal{U}$ from (i) is a finite flat morphism.

The following lemma will be used in the proof of Theorems 4.2 and 5.4.

**3.3. Lemma.** *Let $\mathcal{U}$ be a regular semi-local scheme of geometric type over an infinite field $k$. Let $(p : \mathcal{X} \to \mathcal{U}, \delta : \mathcal{U} \to \mathcal{X}, \mathfrak{f} \in \Gamma(\mathcal{X}, \mathcal{O}_\mathcal{X}))$ be a perfect triple over $\mathcal{U}$. Let $F : Cp(\mathcal{U}) \to Ab$ be a homotopy invariant functor endowed with transfers which satisfies vanishing property (see 2.4, 2.6 and 2.5). Then for each closed subset $\mathcal{Z}$ of the vanishing locus of $\mathfrak{f}$ the following composite vanishes*
$$F_\mathcal{Z}(\mathcal{X}) \xrightarrow{F(\delta)} F_{\delta^{-1}(\mathcal{Z})}(\mathcal{U}) \xrightarrow{F(\mathrm{id}_\mathcal{U})} F_{p(\mathcal{Z})}(\mathcal{U})$$

**3.4. Remark.** The mentioned composite is the map induced by the morphism $\delta : (\mathcal{U}, p(\mathcal{Z})) \to (\mathcal{X}, \mathcal{Z})$ in the category $Cp(\mathcal{U})$. Observe that we have $\delta^{-1}(\mathcal{Z}) \subset p(\mathcal{Z})$, where $p(\mathcal{Z})$ is closed by (ii) of 3.1.

*Proof.* Consider the commutative diagram in the category $Cp(\mathcal{U})$

$$\begin{array}{ccc} (\mathcal{X}, \mathcal{Y}) & \xrightarrow{\mathrm{id}_\mathcal{X}} & (\mathcal{X}, \mathcal{Z}) \\ \delta \uparrow & & \uparrow \delta \\ (\mathcal{U}, Z') & \xrightarrow{\mathrm{id}_\mathcal{U}} & (\mathcal{U}, Z) \end{array}$$

where $Z = \delta^{-1}(\mathcal{Z})$, $Z' = p(\mathcal{Z})$ and $\mathcal{Y} = p^{-1}(Z')$. It gives the relation $F(\mathrm{id}_\mathcal{U}) \circ F(\delta) = F(\delta) \circ F(\mathrm{id}_\mathcal{X})$. Thus to prove the theorem it suffices to check that the following composite vanishes

$$F_\mathcal{Z}(\mathcal{X}) \xrightarrow{F(\mathrm{id}_\mathcal{X})} F_\mathcal{Y}(\mathcal{X}) \xrightarrow{F(\delta)} F_{Z'}(\mathcal{U}) \qquad (\dagger)$$

By Lemma 3.5 below applied to the perfect triple $(\mathcal{X}, \delta, \mathfrak{f})$ we can choose the finite surjective morphism $\pi : \mathcal{X} \to \mathbb{A}^1 \times \mathcal{U}$ from (i) of 3.1 in such a way that it's fibres at the points 0 and 1 of $\mathbb{A}^1$ look as follows:

(a) $\pi^{-1}(\{0\} \times \mathcal{U}) = \delta(\mathcal{U}) \amalg \mathcal{D}_0$ (scheme-theoretically) and $\mathcal{D}_0 \subset \mathcal{X}_\mathfrak{f}$;

(b) $\pi^{-1}(\{1\} \times \mathcal{U}) = \mathcal{D}_1$ and $\mathcal{D}_1 \subset \mathcal{X}_\mathfrak{f}$.

Observe that $\mathcal{Y} = \pi^{-1}(\mathbb{A}^1 \times Z')$. Let $\mathcal{Z}'_0 = \pi^{-1}(\{0\} \times Z') \cap \mathcal{D}_0$ and $\mathcal{Z}'_1 = \pi^{-1}(\{1\} \times Z')$ be the closed subschemes of $\mathcal{Y}$. By definition $\mathcal{Z}'_0$, $\mathcal{Z}'_1$ are the closed subschemes of $\mathcal{D}_0$ and $\mathcal{D}_1$ respectively. Since $\mathcal{Z}$ is contained in the vanishing locus of $\mathfrak{f}$ and $\mathcal{D}_0, \mathcal{D}_1 \subset \mathcal{X}_\mathfrak{f}$ we have $\mathcal{Z} \cap \mathcal{D}_0 = \mathcal{Z} \cap \mathcal{D}_1 = \emptyset$.

The latter means that there are two commutative diagrams in the category $Cp(\mathcal{U})$

$$\begin{array}{ccc} (\mathcal{X}, \mathcal{Y}) & \xrightarrow{\mathrm{id}_\mathcal{X}} & (\mathcal{X}, \mathcal{Z}) \\ I_0 \uparrow & & \uparrow I_0 \\ (\mathcal{D}_0, \mathcal{Z}'_0) & \xrightarrow{\mathrm{id}_{\mathcal{D}_0}} & (\mathcal{D}_0, \emptyset) \end{array} \quad \text{and} \quad \begin{array}{ccc} (\mathcal{X}, \mathcal{Y}) & \xrightarrow{\mathrm{id}_\mathcal{X}} & (\mathcal{X}, \mathcal{Z}) \\ I_1 \uparrow & & \uparrow I_1 \\ (\mathcal{D}_1, \mathcal{Z}'_1) & \xrightarrow{\mathrm{id}_{\mathcal{D}_1}} & (\mathcal{D}_1, \emptyset) \end{array}$$

where $I_0$, $I_1$ are the closed embeddings $\mathcal{D}_0 \hookrightarrow \mathcal{X}$ and $\mathcal{D}_1 \hookrightarrow \mathcal{X}$ respectively.

By vanishing property 2.5 we have $F(\mathcal{D}_0, \emptyset) = F(\mathcal{D}_1, \emptyset) = 0$. Then applying $F$ to the diagrams we immediately get

$$F(I_0) \circ F(\mathrm{id}_\mathcal{X}) = 0 \text{ and } F(I_1) \circ F(\mathrm{id}_\mathcal{X}) = 0. \tag{1}$$

Let $i_0, i_1 : \mathcal{U} \hookrightarrow \mathbb{A}^1 \times \mathcal{U}$ be the closed embeddings which correspond to the points 0 and 1 of $\mathbb{A}^1$ respectively. The homotopy invariance property 2.4 implies that

$$F(i_0) = F(i_1) : F_{\mathbb{A}^1 \times Z'}(\mathbb{A}^1 \times \mathcal{U}) \to F_{Z'}(\mathcal{U}). \tag{2}$$

The base change property 2.6.(i) applied to the fibred product diagram

$$\begin{array}{ccc} (\mathcal{X}, \mathcal{Y}) & \xleftarrow{I_1} & (\mathcal{D}_1, \mathcal{Z}'_1) \\ \pi \downarrow & & \downarrow \pi \\ (\mathbb{A}^1 \times \mathcal{U}, \mathbb{A}^1 \times Z') & \xleftarrow{i_1} & (\mathcal{U}, Z') \end{array}$$

gives the relation

$$F(i_1) \circ \mathrm{Tr}^{\mathcal{X}}_{\mathbb{A}^1 \times \mathcal{U}} = \mathrm{Tr}^{\mathcal{D}_1}_{\mathcal{U}} \circ F(I_1), \tag{3}$$

where $\mathrm{Tr}^{\mathcal{X}}_{\mathbb{A}^1 \times \mathcal{U}} : F_\mathcal{Y}(\mathcal{X}) \to F_{\mathbb{A}^1 \times Z'}(\mathbb{A}^1 \times \mathcal{U})$ and $\mathrm{Tr}^{\mathcal{D}_1}_{\mathcal{U}} : F_{\mathcal{Z}'_1}(\mathcal{D}_1) \to F_{Z'}(\mathcal{U})$ are the transfer maps for the finite flat morphism $\pi$ and $\pi|_{\mathcal{D}_1}$ respectively.

Consider the commutative diagram

$$\begin{array}{ccccccc} F_\mathcal{Z}(\mathcal{X}) & \xrightarrow{F(\mathrm{id}_\mathcal{X})} & F_\mathcal{Y}(\mathcal{X}) & \xrightarrow{(F(\delta), F(I_0))} & F_{Z' \amalg \mathcal{Z}'_0}(\mathcal{U} \amalg \mathcal{D}_0) & \xrightarrow{\text{``+''}} & F_{Z'}(\mathcal{U}) \oplus F_{\mathcal{Z}'_0}(\mathcal{D}_0) \\ & & \mathrm{Tr}^{\mathcal{X}}_{\mathbb{A}^1 \times \mathcal{U}} \downarrow & & \mathrm{Tr} \downarrow & & \swarrow \mathrm{id} + \mathrm{Tr}^{\mathcal{D}_0}_{\mathcal{U}} \\ & & F_{\mathbb{A}^1 \times Z'}(\mathbb{A}^1 \times \mathcal{U}) & \xrightarrow{F(i_0)} & F_{Z'}(\mathcal{U}) & & \end{array}$$

$$\tag{4}$$

where the central square commutes by 2.6.(i) and the right triangle commutes by 2.6.(ii). In the diagram we identify $\mathcal{U}$ with $\delta(\mathcal{U})$ by means of the isomorphism $\delta : \mathcal{U} \to \delta(U)$ and use the property 2.6.(iii) to identify $\mathrm{Tr}_{\mathcal{U}}^{\delta(\mathcal{U})}$ with $F(\delta)$.

The following chain of relations shows that the composite (†) vanishes and we finish the proof of the lemma.

$$F(\delta) \circ F(\mathrm{id}_{\mathcal{X}}) \stackrel{(1)}{=} (\mathrm{id} + \mathrm{Tr}_{\mathcal{U}}^{\mathcal{D}_0}) \circ (F(\delta), F(I_0)) \circ F(\mathrm{id}_{\mathcal{X}}) \stackrel{(4)}{=} F(i_0) \circ \mathrm{Tr}_{\mathbb{A}^1 \times \mathcal{U}}^{\mathcal{X}} \circ F(\mathrm{id}_{\mathcal{X}}) \stackrel{(2)}{=}$$

$$F(i_1) \circ \mathrm{Tr}_{\mathbb{A}^1 \times \mathcal{U}}^{\mathcal{X}} \circ F(\mathrm{id}_{\mathcal{X}}) \stackrel{(3)}{=} \mathrm{Tr}_{\mathcal{U}}^{\mathcal{D}_1} \circ F(I_1) \circ F(\mathrm{id}_{\mathcal{X}}) \stackrel{(1)}{=} 0$$

□

The following lemma is the semi-local version of Geometric Presentation Lemma ([8], 10.1)

**3.5. Lemma.** *Let $R$ be a semi-local essentially smooth algebra over an infinite field $k$ and $A$ an essentially smooth $k$-algebra, which is finite over the polynomial algebra $R[t]$. Suppose that $e : A \to R$ is an $R$-augmentation and let $I = \ker e$. Assume that $A$ is smooth over $R$ at every prime containing $I$. Given $f \in A$ such that $A/Af$ is finite over $R$ we can find an $s \in A$ such that*

1. *$A$ is finite over $R[s]$.*

2. *$A/As = A/I \times A/J$ for some ideal $J$ of $A$.*

3. *$J + Af = A$.*

4. *A(s-1)+Af=A.*

*Proof.* In the proof of ([8], 10.1) replace the reduction modulo maximal ideal by the reduction modulo radical of the semi-local ring. □

# 4 The constant case of the Effacement Theorem

The following important definition comes from ([3], 2.1.1)

**4.1. Definition.** Let $X$ be a smooth affine variety over a field $k$. Let $x = \{x_1, \ldots, x_n\}$ be a finite set of points of $X$ and let $\mathcal{U} = \operatorname{Spec} \mathcal{O}_{X,x}$ be the semi-local scheme at $x$. A contravariant functor $F : Cp(X) \to Ab$ is effaceable at $x$ if the following condition satisfied:

Given $m \geq 1$, for any closed subset $Z \subset X$ of codimension $m$, there exist a closed subset $Z' \subset \mathcal{U}$ such that

(1) $Z' \supset Z \cap \mathcal{U}$ and $\operatorname{codim}_\mathcal{U}(Z') \geq m - 1$;

(2) the composite $F_Z(X) \xrightarrow{F(j)} F_{Z \cap \mathcal{U}}(\mathcal{U}) \xrightarrow{F(\operatorname{id}_\mathcal{U})} F_{Z'}(\mathcal{U})$ vanishes, where $j : \mathcal{U} \to X$ is the canonical embedding and $Z \cap \mathcal{U} = j^{-1}(Z)$.

**4.2. Theorem.** *Let $X$ be a smooth affine variety over an infinite field $k$ and $x \subset X$ be a finite set of points. Let $G : Cp(k) \to Ab$ be a homotopy invariant functor endowed with transfers which satisfies vanishing property (see 2.4, 2.6 and 2.5). Let $F = p^*G$ denote the restriction of $G$ to $Cp(X)$ by means of the structural morphisms $p : X \to \operatorname{Spec} k$. Then $F$ is effaceable at $x$.*

**4.3. Remark.** As a consequence of this theorem we immediately get Theorem 1.1 for the case when the complex $\mathcal{C}$ comes from the base field $k$, i.e., each sheaf in $\mathcal{C}$ can be represented as $p^*G$ for some sheaf $G$ on $(\operatorname{Spec} k)_{et}$, where $p : X \to \operatorname{Spec} k$ is the structural morphism.

Indeed, the étale cohomology functor $F(X, Z) = H_Z^*(X, \mathcal{C})$ satisfies all the hypotheses of Theorem 4.2. It is homotopy invariant according to ([3], 7.3.(1)). It satisfies vanishing property by the very definition and has transfer maps by section 6. So that we may apply Theorem 4.2 to the functor $F$ and get $F$ is effaceable. As it was mentioned before Theorem 1.1 now follows from Lemma 5.7.

*Proof.* Here we give only the proof for the case of closed points. The general case can be reduced to that one as in the proof of Theorem 7.1, [9].

We may assume $x \cap Z$ in non-empty. Indeed, if $x \cap Z = \emptyset$ then the theorem follows by the vanishing property of $F$ (see 2.5).

Let $f \neq 0$ be a regular function on $X$ such that $Z$ is a closed subset of the vanishing locus of $f$. By Quillen's trick ([12], 5.12), ([11], 1.2) we can find a morphism $q : X \to \mathbb{A}^{n-1}$, where $n = \dim X$, such that

(a) $q|_{f=0} : \{f = 0\} \to \mathbb{A}^{n-1}$ is a finite morphism;

(b) $q$ is smooth at the points $x$;

(c) $q$ can be factorized as $q = pr \circ \Pi$, where $\Pi : X \to \mathbb{A}^n$ is a finite surjective morphism and $pr : \mathbb{A}^n \to \mathbb{A}^{n-1}$ is a linear projection.

Consider the base change diagram for the morphism $q$ by means of the composite $r : \mathcal{U} = \operatorname{Spec} \mathcal{O}_{X,x} \xrightarrow{j} X \xrightarrow{q} \mathbb{A}^{n-1}$.

$$\begin{array}{ccc} \mathcal{X} & \xrightarrow{r_X} & X \\ p \downarrow & & \downarrow q \\ \mathcal{U} & \xrightarrow{r} & \mathbb{A}^{n-1} \end{array}$$

So we have $\mathcal{X} = \mathcal{U} \times_{\mathbb{A}^{n-1}} X$ and $p$, $r_X$ denote the canonical projections on $\mathcal{U}$, $X$ respectively. Let $\delta : \mathcal{U} \to \mathcal{X} = \mathcal{U} \times_{\mathbb{A}} X$ be the diagonal embedding. Clearly $\delta$ is a section of $p$. Set $\mathfrak{f} = r_X^*(f)$. Take instead of $\mathcal{X}$ it's irreducible component containing $\delta(\mathcal{U})$ and instead of $\mathfrak{f}$ it's restriction to this irreducible component (since $x \cap Z$ is non-empty the vanishing locus of $\mathfrak{f}$ on the component containing $\delta(\mathcal{U})$ is non-empty as well).

Now assuming the triple $(p : \mathcal{X} \to \mathcal{U}, \delta, \mathfrak{f})$ is a perfect triple over $\mathcal{U}$ (see 3.1 for the definition) we complete the proof as follows:

Let $\mathcal{Z} = r_X^{-1}(Z)$ be a closed subset of the vanishing locus of $\mathfrak{f}$. Let $Z' = p(\mathcal{Z})$ be a closed subscheme of $\mathcal{U}$. Since $r_X \circ \delta = j$ we have $\delta^{-1}(\mathcal{Z}) = j^{-1}(Z) = Z \cap \mathcal{U}$. By Specialization Lemma 3.3 applied to the perfect triple $(\mathcal{X}, \delta, \mathfrak{f})$ and the functor $j^*F : Cp(\mathcal{U}) \to Ab$ the composite

$$F_{\mathcal{Z}}(\mathcal{X}) \xrightarrow{F(\delta)} F_{Z \cap \mathcal{U}}(\mathcal{U}) \xrightarrow{F(\operatorname{id}_{\mathcal{U}})} F_{Z'}(\mathcal{U})$$

vanishes. In particular, the composite

$$F_Z(X) \xrightarrow{F(j)} F_{Z \cap \mathcal{U}}(\mathcal{U}) \xrightarrow{F(\operatorname{id}_{\mathcal{U}})} F_{Z'}(\mathcal{U}) \qquad (*)$$

vanishes as well. Clearly $Z' \supset Z \cap \mathcal{U}$. By 3.1.(i) we have $\dim \mathcal{X} = \dim \mathcal{U} + 1$. On the other hand the morphism $r_\mathcal{X} : \mathcal{X} \to X$ is flat (even essentially smooth) and, thus, $\operatorname{codim}_\mathcal{X}(\mathcal{Z}) = \operatorname{codim}_X Z$. Therefore, we have $\operatorname{codim}_\mathcal{U}(Z') = \operatorname{codim}_\mathcal{X}(\mathcal{Z}) - 1 = m - 1$. □

Hence, it remains to prove the following:

**4.4. Lemma.** *The triple $(p : \mathcal{X} \to \mathcal{U}, \delta, \mathfrak{f})$ is perfect over $\mathcal{U}$.*

*Proof.* By the property (c) one has $q = pr \circ \Pi$ with a finite surjective morphism $\Pi : X \to \mathbb{A}^n$ and a linear projection $pr : \mathbb{A}^n \to \mathbb{A}^{n-1}$. Taking the base change of $\Pi$ by means of $r : \mathcal{U} \to \mathbb{A}^{n-1}$ one gets a finite surjective $\mathcal{U}$-morphism $\pi : \mathcal{X} \to \mathbb{A}^1 \times \mathcal{U}$. This checks (i) of 3.1. Since the closed subset $\{f = 0\}$ of $X$ is finite over $\mathbb{A}^{n-1}$ the closed subset $\{\mathfrak{f} = 0\}$ of $\mathcal{X}$ is finite over $\mathcal{U}$ and we get 3.1.(ii). Since $q$ is smooth at $x$ the morphism $r : \mathcal{U} \to \mathbb{A}^{n-1}$ is essentially smooth. Thus the morphism $r_X : \mathcal{X} \to X$ is essentially smooth as the base change of the morphism $r$. The variety $X$ is smooth over $k$ implies that $\mathcal{X}$ is essentially smooth over $k$ as well. Since $q$ is smooth at $x$ the morphism $p$ is smooth at each point $y \in \mathcal{X}$ with $r_X(y) \in x$. In particular $p$ is smooth at the points $\delta(x_i)$ ($x_i \in \mathcal{U}$). Since $\mathcal{U}$ is semi-local $\delta(\mathcal{U})$ is semi-local and $p$ is smooth along $\delta(\mathcal{U})$. This checks (iii) of 3.1. Since $\mathcal{X}$ is irreducible 3.1.(iv) holds. And we have proved the lemma and the theorem. □

## 5 The Effacement Theorem

To prove Theorem 1.1 in the general case we have to put additional condition on our functor $F$. In order to formulate this condition we need some notation.

**5.1. Notation.** Let $\rho : Y \to X \times X$ be a finite étale morphism together with a section $s : X \to Y$ over the diagonal embedding $\Delta : X \to X \times X$, i.e., $\rho \circ s = \Delta$. Let $pr_1, pr_2 : X \times X \to X$ be the canonical projections. We denote $p_1, p_2 : Y \to X$ to be the composite $pr_1 \circ \rho$, $pr_2 \circ \rho$ respectively.

For a contravariant functor $F : Cp(X) \to Ab$ consider it's pull-backs $p_1^*F$ and $p_2^*F : Cp(Y) \to Ab$ by means of $p_1$ and $p_2$ respectively. From this point on we denote $F_1 = p_1^*F$ and $F_2 = p_2^*F$. By definition we have

$$F_i(Y' \to Y, Z) = F(Y' \to Y \xrightarrow{p_i} X, Z).$$

**5.2. Remark.** In general case the functors $F_1$ and $F_2$ are not equivalent. Moreover, the functors $pr_1^*F$ and $pr_2^*F$ are different. But in the case when $F$ comes from the base field $k$, i.e., $F = p^*G$ where $p : X \to \operatorname{Spec} k$ is the structural morphism and $G : Cp(k) \to Ab$ is a contravariant functor, these functors coincide with each other.

**5.3. Definition.** We say a contravariant functor $F : Cp(X) \to Ab$ has a finite monodromy of the type $(\rho : Y \to X \times X, s : X \to Y)$, where $\rho$ is a finite étale morphism and $s$ is a section of $\rho$ over the diagonal, if there exists an isomorphism $\Phi : F_1 \to F_2$ of functors on $Cp(Y)$. A functor $F : Cp(X) \to Ab$ is said to be a functor with finite monodromy if $F$ has a finite monodromy of some type.

**5.4. Theorem.** *Let $X$ be a smooth affine variety over an infinite field $k$ and $x \subset X$ be a finite set of points. Let $F : Cp(X) \to Ab$ be a homotopy invariant functor endowed with transfers which satisfies vanishing property (see 2.4, 2.6 and 2.5). If $F$ is a functor with finite monodromy then $F$ is effaceable at $x$.*

*Proof.* Similar to the proof of Theorem 4.2 let $f \neq 0$ be a regular function on $X$ such that $Z$ is a closed subset of the vanishing locus of $f$. We may assume $x \cap Z$ is non-empty. Consider the fibred product diagram from the proof of Theorem 4.2

$$\begin{array}{ccc} \mathcal{X} & \xrightarrow{r_X} & X \\ p \downarrow & & \downarrow q \\ \mathcal{U} & \xrightarrow{r} & \mathbb{A}^{n-1} \end{array}$$

We have the projection $p : \mathcal{X} = \mathcal{U} \times_{\mathbb{A}^{n-1}} X \to \mathcal{U}$, the section $\delta : \mathcal{U} \to \mathcal{X}$ of $p$ and the regular function $\mathfrak{f} = r_X^*(f)$.

Since $F$ is the functor with finite monodromy there is a finite étale morphism $\rho : Y \to X \times X$, a section $s : X \to Y$ of $\rho$ over the diagonal embedding and a functor isomorphism $\Phi : F_1 \to F_2$ as in 5.1 and 5.3. Consider the base change diagram for the morphism $\rho : Y \to X \times X$ by means of the composite $g : \mathcal{X} \xrightarrow{(p, r_X)} \mathcal{U} \times X \xrightarrow{(j, \mathrm{id})} X \times X$

$$\begin{array}{ccc} \tilde{\mathcal{X}} & \xrightarrow{\tilde{g}} & Y \\ \tilde{\rho} \downarrow & & \downarrow \rho \\ \mathcal{X} & \xrightarrow{g} & X \times X \end{array}$$

Then $\tilde{\rho}$ is a finite étale morphism and there is the section $\tilde{\delta} : \mathcal{U} \to \tilde{\mathcal{X}}$ of the composite $\tilde{p} = p \circ \tilde{\rho} : \tilde{\mathcal{X}} \to \mathcal{U}$ such that $\tilde{\rho} \circ \tilde{\delta} = \delta$ ($\tilde{\delta}$ is the base change of the morphism $s : X \to Y$ by means of $\tilde{g} : \tilde{\mathcal{X}} \to Y$). Set $\tilde{\mathfrak{f}} = \tilde{\rho}^{-1}(\mathfrak{f})$. As in the proof of 4.2 we replace $\tilde{\mathcal{X}}$ by it's irreducible component containing $\tilde{\delta}(\mathcal{U})$ and $\tilde{\mathfrak{f}}$ by it's restriction to this component. By Lemma 5.5 below the triple $(\tilde{p} : \tilde{\mathcal{X}} \to \mathcal{U}, \tilde{\delta}, \tilde{\mathfrak{f}})$ is perfect.

Let $\mathcal{Z} = r_X^{-1}(Z)$ and $\tilde{\mathcal{Z}} = \tilde{\rho}^{-1}(\mathcal{Z})$. Set $Z' = p(\mathcal{Z})$. Observe that $\tilde{p}(\tilde{\mathcal{Z}}) = Z'$. The commutative diagram

$$\begin{array}{ccccc} F_Z(X) & \xrightarrow{F(j)} & F_{Z\cap\mathcal{U}}(\mathcal{U}) & \xrightarrow{F(\mathrm{id}_\mathcal{U})} & F_{Z'}(\mathcal{U}) \\ {\scriptstyle F(r_X)}\downarrow & {\scriptstyle F(\delta)}\nearrow & & {\scriptstyle F(\tilde{\delta})}\uparrow & \\ F_\mathcal{Z}(\mathcal{X}) & \xrightarrow{F(\tilde{\rho})} & F_{\tilde{\mathcal{Z}}}(\tilde{\mathcal{X}}) & & \end{array}$$

shows that to prove the relation $F(\mathrm{id}_\mathcal{U}) \circ F(j) = 0$ (compare with $(*)$ of the proof of 4.2) it suffices to check the relation $F(\mathrm{id}_\mathcal{U}) \circ F(\tilde{\delta}) = 0$.

Consider the pull-backs of the functors $F_1$, $F_2$ and the functor isomorphism $\Phi$ by means of the morphism $\tilde{g} : \tilde{\mathcal{X}} \to Y$. We shall use the same notation $F_1$, $F_2$ and $\Phi$ for these pull-backs till the end of this proof. So we have $F_1 = \tilde{g}^*(p_1^*F)$ and $F_2 = \tilde{g}^*(p_2^*F)$. The functor isomorphism $\Phi : F_1 \xrightarrow{\cong} F_2$ provides us with the commutative diagram of functors values on $Cp(\tilde{\mathcal{X}})$

$$\begin{array}{ccccc} (F_1)_{\tilde{\mathcal{Z}}}(\tilde{\mathcal{X}}) & \xrightarrow{F_1(\tilde{\delta})} & (F_1)_{Z\cap\mathcal{U}}(\mathcal{U}) & \xrightarrow{F_1(\mathrm{id}_\mathcal{U})} & (F_1)_{Z'}(\mathcal{U}) \\ {\scriptstyle \Phi}\downarrow{\scriptstyle \cong} & & {\scriptstyle \Phi}\downarrow{\scriptstyle \cong} & & {\scriptstyle \cong}\downarrow{\scriptstyle \Phi} \\ (F_2)_{\tilde{\mathcal{Z}}}(\tilde{\mathcal{X}}) & \xrightarrow{F_2(\tilde{\delta})} & (F_2)_{Z\cap\mathcal{U}}(\mathcal{U}) & \xrightarrow{F_2(\mathrm{id}_\mathcal{U})} & (F_2)_{Z'}(\mathcal{U}) \end{array}$$

where the structure of an $\tilde{\mathcal{X}}$-scheme on $\mathcal{U}$ is given by $\tilde{\delta}$.

Since $r_X \circ \tilde{\rho} = p_2 \circ \tilde{g}$ we have $\tilde{\rho}^*(r_X^*F) = F_2$. Thus to check the relation $F(\mathrm{id}_\mathcal{U}) \circ F(\tilde{\delta}) = 0$ for the functor $F$ we have to verify the same relation $F_2(\mathrm{id}_\mathcal{U}) \circ F_2(\tilde{\delta}) = 0$ for the functor $F_2$. Then by commutativity of the diagram it suffices to prove the relation $F_1(\mathrm{id}_\mathcal{U}) \circ F_1(\tilde{\delta}) = 0$ for the functor $F_1$.

Since $j \circ \tilde{p} = p_1 \circ \tilde{g}$ we have $\tilde{p}^*(j^*F) = F_1 : Cp(\tilde{\mathcal{X}}) \to Ab$. Thereby it suffices to prove the relation $G(\mathrm{id}_\mathcal{U}) \circ G(\tilde{\delta}) = 0$ for the functor $G = j^*F : Cp(\mathcal{U}) \to Ab$. This relation follows immediately from Theorem 3.3 applied to the functor $G$, the triple $(\tilde{p} : \tilde{\mathcal{X}} \to \mathcal{U}, \tilde{\delta}, \tilde{\mathfrak{f}})$ and the closed subset $\tilde{\mathcal{Z}} \subset \tilde{\mathcal{X}}$. $\square$

**5.5. Lemma.** *The triple $(\tilde{\mathcal{X}}, \tilde{\delta}, \tilde{\mathfrak{f}})$ is perfect over $\mathcal{U}$.*

*Proof.* Observe that the triple $(p : \mathcal{X} \to \mathcal{U}, \delta, \mathfrak{f})$ is perfect by Lemma 4.4, the morphism $\tilde{\rho} : \tilde{\mathcal{X}} \to \mathcal{X}$ is finite étale and $\tilde{\rho} \circ \tilde{\delta} = \delta$ for the section $\tilde{\delta} : \mathcal{U} \to \tilde{\mathcal{X}}$ of the morphism $\tilde{p} : \tilde{\mathcal{X}} \to \mathcal{U}$. For the finite surjective morphism of $\mathcal{U}$-schemes $\pi : \mathcal{X} \to \mathbb{A}^1 \times \mathcal{U}$ the composite $\tilde{\mathcal{X}} \xrightarrow{\tilde{\rho}} \mathcal{X} \xrightarrow{\pi} \mathbb{A}^1 \times \mathcal{U}$ is a finite surjective morphism of $\mathcal{U}$-schemes as well. This proves 3.1.(i). Since $\tilde{\rho}$ is finite and the vanishing locus of $\mathfrak{f}$ is finite over $\mathcal{U}$ the vanishing locus of the function $\tilde{\mathfrak{f}}$ is finite over $\mathcal{U}$ as well. This proves 3.1.(ii). Since $\tilde{\rho} \circ \tilde{\delta} = \delta$, $\tilde{\rho}$ is étale and $p$ is smooth along $\delta(\mathcal{U})$ the morphism $\tilde{p}$ is smooth along $\tilde{\delta}(\mathcal{U})$. Since the scheme $\mathcal{X}$ is essentially smooth over $k$ and $\tilde{\rho}$ is étale the scheme $\tilde{\mathcal{X}}$ is essentially smooth over $k$. This proves 3.1.(iii). Since $\tilde{\mathcal{X}}$ is irreducible we have 3.1.(iv). And the lemma is proven. □

**5.6. Remark.** In order to proof Theorem 1.1 we apply the same arguments as in Remark 4.3. The only additional fact we have to check is that the étale cohomology functor $F(X, Z) = H_Z^*(X, \mathcal{C})$ is a functor with finite monodromy. This is done in the section 7 (see Corollary 7.6).

Similar to Proposition 2.1.2, [3], we have

**5.7. Lemma.** *Let $\mathcal{U}$ be a semi-local regular scheme of geometric type over a field $k$, i.e., $\mathcal{U} = \operatorname{Spec} \mathcal{O}_{X,x}$ for some smooth affine variety $X$ and a finite set of points $x = \{x_1, \ldots, x_n\}$ of $X$. Suppose the étale cohomology functor $F(Y, Z) = H_Z^*(Y, \mathcal{C})$ from Theorem 1.1 is effaceable at $x$. Then, in the exact couple ([3], 1.1) defining the coniveau spectral sequence for $(\mathcal{U}, \mathcal{C})$, the map $i^{p,q}$ is identically $0$ for all $p > 0$. In particular, we have $E_2^{p,q} = H^q(\mathcal{U}, \mathcal{C})$ if $p = 0$ and $E_2^{p,q} = 0$ if $p > 0$. And the Cousin complex ([3], 1.3) yields the exact sequence from Theorem 1.1.*

*Proof.* Consider the commutative diagram

$$\begin{array}{ccccc} H_Z^n(X, \mathcal{C}) & \longrightarrow & H_{Z \cap \mathcal{U}}^n(\mathcal{U}, \mathcal{C}) & \longrightarrow & H_{Z'}^n(\mathcal{U}, \mathcal{C}) \\ \downarrow & & \downarrow & & \downarrow \\ H_{X^{(m)}}^n(X, \mathcal{C}) & \longrightarrow & H_{\mathcal{U}^{(m)}}^n(\mathcal{U}, \mathcal{C}) & \longrightarrow & H_{\mathcal{U}^{(m-1)}}^n(\mathcal{U}, \mathcal{C}) \end{array}$$

The composition of arrows in the first row is identically $0$ for any $n$. Therefore the compositions $H_Z^n(X, \mathcal{C}) \to H_{\mathcal{U}^{(m)}}^n(\mathcal{U}, \mathcal{C}) \to H_{\mathcal{U}^{(m-1)}}^n(\mathcal{U}, \mathcal{C})$ are $0$. Passing to the limit over $Z$, this gives that the compositions $H_{X^{(m)}}^n(X, \mathcal{C}) \to H_{\mathcal{U}^{(m)}}^n(\mathcal{U}, \mathcal{C}) \to H_{\mathcal{U}^{(m-1)}}^n(\mathcal{U}, \mathcal{C})$ are $0$. Passing to the limit over open neighborhoods of $x$, we get that the map $i^{m,n-m} : H_{\mathcal{U}^{(m)}}^n(\mathcal{U}, \mathcal{C}) \to H_{\mathcal{U}^{(m-1)}}^n(\mathcal{U}, \mathcal{C})$ is itself $0$ for any $m \geq 1$. □

# 6 Transfers for Étale Cohomology

**6.1. Notation.** Let $\pi : Y \to X$ be a finite flat morphism of schemes and let $F$ be a sheaf on the big étale site $Et/X$. For an $X$-scheme $X'$ set $Y' = X' \times_X Y$ and denote the projection $Y' \to X'$ by $\pi'$. If $Y' = Y'_1 \amalg Y'_2$ (disjoint union) then set $\pi'_i = \pi'|_{Y'_i}$. If $X'' \xrightarrow{g} X'$ is an $X$-scheme morphism then set $Y'' = X'' \times_X Y$ and denote by $\pi'' : Y'' \to X''$ the projection on $X''$. Denote by $g_Y : Y'' \to Y'$ the morphism $g \times \mathrm{id}_Y$. If $Y' = Y'_1 \amalg Y'_2$ then set $Y''_i = g_Y^{-1}(Y'_i)$ and define $g_{Y,i} : Y''_i \to Y'_i$ as the restriction of $g_Y$.

If $Z \subset X$ is a closed subset then set $S = \pi^{-1}(Z)$, $Z' = X' \times_X Z$, $Z'' = X'' \times_X Z$, $S' = (\pi')^{-1}(Z')$, $S'' = (\pi'')^{-1}(Z'')$. If $Y' = Y'_1 \amalg Y'_2$ then set $S'_i = Y'_i \cap S'$, $S''_i = Y''_i \cap S''$.

**6.2. Trace Maps.** Deligne in [4] constructed trace maps for finite flat morphisms. In particular for each $X$-scheme $X'$ and for every presentation of the scheme $Y'$ in the form $Y' = Y'_1 \amalg Y'_2$ there are certain trace maps $\mathrm{Tr}_{\pi'_i} : \Gamma(Y'_i, F) \to \Gamma(X', F)$. These maps satisfy the following properties

(i) (base change) the diagram

$$\begin{array}{ccc} \Gamma(Y''_i, F) & \xleftarrow{g^*_{Y,i}} & \Gamma(Y'_i, F) \\ {\scriptstyle \mathrm{Tr}_{\pi''_i}} \downarrow & & \downarrow {\scriptstyle \mathrm{Tr}_{\pi'_i}} \\ \Gamma(X'', F) & \xleftarrow{g^*_Y} & \Gamma(X', F) \end{array}$$

commutes;

(ii) (additivity) the diagram

$$\begin{array}{ccc} \Gamma(Y', F) & \xrightarrow{\text{``}+\text{''}} & \Gamma(Y'_1, F) \oplus \Gamma(Y'_2, F) \\ {\scriptstyle \mathrm{Tr}_{\pi'}} \downarrow & & \downarrow {\scriptstyle \mathrm{Tr}_{\pi'_1} + \mathrm{Tr}_{\pi'_2}} \\ \Gamma(X', F) & \xrightarrow{\mathrm{id}} & \Gamma(X', F) \end{array}$$

commutes;

(iii) (normalisation) if $\pi'_1 : Y'_1 \to X'$ is an isomorphism then the composite map

$$\Gamma(X', F) \xrightarrow{(\pi'_1)^*} \Gamma(Y'_1, F) \xrightarrow{\mathrm{Tr}_{\pi'_1}} \Gamma(X', F)$$

is the identity;

(iv) maps $\text{Tr}_{\pi'_i}$ are functorial with respect to sheaves $F$ on $Et/X$.

**6.3. Trace Maps with Supports.** If $Z \subset X$ is a closed subset then for an $X$-scheme $X'$ we denote $\Gamma_{Z'}(X', F) = \ker(\Gamma(X', F) \to \Gamma(X' - Z', F))$ and if $Y' = Y'_1 \amalg Y'_2$ we denote $\Gamma_{S'_i}(Y'_i, F) = \ker(\Gamma(Y'_i, F) \to \Gamma(Y'_i - S'_i, F))$. Thus the trace maps induce certain maps $\text{Tr}_{\pi'_i} : \Gamma_{S'_i}(Y'_i, F) \to \Gamma_{Z'}(X', F)$. These maps clearly satisfy the following properties

(i) (base change) the diagram

$$\begin{array}{ccc} \Gamma_{S''_i}(Y''_i, F) & \xleftarrow{g^*_{Y,i}} & \Gamma_{S'_i}(Y'_i, F) \\ \text{Tr}_{\pi''_i} \downarrow & & \downarrow \text{Tr}_{\pi'_i} \\ \Gamma_{Z''}(X'', F) & \xleftarrow{g^*_Y} & \Gamma_{Z'}(X', F) \end{array}$$

commutes;

(ii) (additivity) the diagram

$$\begin{array}{ccc} \Gamma_{S'}(Y', F) & \xrightarrow{\text{``+''}} & \Gamma_{S'_1}(Y'_1, F) \oplus \Gamma_{S'_2}(Y'_2, F) \\ \text{Tr}_{\pi'} \downarrow & & \downarrow \text{Tr}_{\pi'_1} + \text{Tr}_{\pi'_2} \\ \Gamma_{Z'}(X', F) & \xrightarrow{\text{id}} & \Gamma_{Z'}(X', F) \end{array}$$

commutes;

(iii) (normalisation) if $\pi'_1 : Y'_1 \to X'$ is an isomorphism then the composite map
$$\Gamma_{Z'}(X', F) \xrightarrow{(\pi'_1)^*} \Gamma_{S'_1}(Y'_1, F) \xrightarrow{\text{Tr}_{\pi'_1}} \Gamma_{Z'}(X', F)$$
is the identity;

(iv) maps $\text{Tr}_{\pi'_i}$ are functorial with respect to sheaves $F$ on $Et/X$.

**6.4. Transfers for Cohomology.** Let $0 \to F \to \mathfrak{I}^\bullet$ be an injective resolution of the sheaf $F$ on $Et/X$. Then for a closed subset $Z \subset X$ and for a presentation $Y' = Y'_1 \amalg Y'_2$ one has

$$H^p_{S'_i}(Y'_i, F) := H^p(\Gamma_{S'_i}(Y'_i, \mathfrak{I}^\bullet)) \, , \, H^p_{Z'}(X', F) := H^p(\Gamma_{Z'}(X', \mathfrak{I}^\bullet)).$$

Thereby 6.3.(iv) shows that the trace maps $\mathrm{Tr}_{\pi_i'} : \Gamma_{S_i'}(Y_i', \mathfrak{I}^r) \to \Gamma_{Z'}(X', \mathfrak{I}^r)$ determine a morphism of complexes $\Gamma_{S_i'}(Y_i', \mathfrak{I}^\bullet) \to \Gamma_{Z'}(X', \mathfrak{I}^\bullet)$. Thus one gets the induced map which we will denote by $H^p(\mathrm{Tr}_{\pi_i'}) : H^p_{S_i'}(Y_i', F) \to H^p_{Z'}(X', F)$. These trace maps satisfy the following properties (similar to the properties 6.3)

(i) the base changing property;

(ii) the additivity property;

(iii) the normalisation property;

(iv) the functorality with respect to sheaves on $Et/X$.

Thus, the following lemma is proven

**6.5. Lemma.** *For a sheaf $F$ on the small étale site $X_{et}$ the functor $(X', Z') \mapsto H^p_{Z'}(X', F)$ is a functor $Cp(X) \to Ab$ endowed with transfers in the sense of 2.6.*

# 7 The Étale Cohomology functor is a functor with finite monodromy

The goal of this section is to check that the étale cohomology functor with supports $(X, Z) \to H^*_Z(X, \mathcal{C})$ from Theorem 1.1 is a functor with finite monodromy (see Definition 5.3). In order to do this, first, we prove Lemma 7.3. The desired fact (Corollary 7.6) then follows immediately.

**7.1. Notation.** Let $k$ be a field and let $X$ be a smooth affine variety over $k$. Let $\tilde{X} \to X$ be an étale Galois covering with the Galois group $G$. We will assume in this section that $\tilde{X}$ is irreducible. Denote by $C(\tilde{X}/X)$ a category consisting of those locally constant constructible sheaves on the small étale site $X_{et}$ which become constant sheaves over $\tilde{X}_{et}$. The morphisms $p_i$ induce functors $p_i^*$ from the category $C(\tilde{X}/X)$ to the category of sheaves on $Y_{et}$.

**7.2. Notation.** To state Lemma 7.3 below we use an explicit construction of a scheme $Y$ (see 5.3) suggested by H. Esnault:

Set $Y = (\tilde{X} \times \tilde{X})/G$, where $G$ acts diagonally on $\tilde{X} \times \tilde{X}$. The canonical projections $pr_i : \tilde{X} \times \tilde{X} \to \tilde{X}$ modulo $G$-action determine the maps $p_i :$

$Y \to X$ and, thus, the morphism $\rho : Y \to X \times X$. The diagonal embedding $\tilde{X} \hookrightarrow \tilde{X} \times \tilde{X}$ modulo $G$-action determines a closed embedding $\delta : X \hookrightarrow Y$. Clearly $\delta$ is a section of the morphism $\rho : Y \to X \times X$ over the diagonal embedding $\Delta : X \to X \times X$.

**7.3. Lemma.** *There is a functor isomorphism $\Phi : p_1^* \xrightarrow{\cong} p_2^*$. In other words there is a family of isomorphisms $\{\Phi_F : p_1^* F \xrightarrow{\cong} p_2^* F \mid F \in C(\tilde{X}/X)\}$ such that for each morphism $F' \xrightarrow{\alpha} F''$ the diagram*

$$\begin{array}{ccc} p_1^* F' & \xrightarrow{p_1^*(\alpha)} & p_1^* F'' \\ \Phi_{F'} \downarrow & & \downarrow \Phi_{F''} \\ p_2^* F' & \xrightarrow{p_2^*(\alpha)} & p_2^* F'' \end{array}$$

*is commutative.*

*Proof.* Let $F$ be a sheaf on $X_{et}$ such that $F|_{\tilde{X}_{et}} = C$ is a constant sheaf. Let $E$ denote the étale group scheme representing $F$ on $X_{et}$ then we may write $E \times_X \tilde{X} = \tilde{X}^C$ and for $W \in X_{et}$ we have

$$\Gamma(W, F) = Mor_X(W, E) = Mor_X^G(W \times_X \tilde{X}, E) = Mor_{\tilde{X}}^G(W \times_X \tilde{X}, \tilde{X}^C). \tag{$*$}$$

Let $V \in Y_{et}$ be an $Y$-scheme. Denote by $V_1$ the $Y$-scheme $V$ considered as an $X$-scheme by means of $V \to Y \xrightarrow{p_1} X$. More precisely, the lower index 1 or 2 in the notation $V_1$ or $V_2$ shows which structure of an $X$-scheme is chosen for a given $Y$-scheme $V$. Then similar to $(*)$ for each $V \in Y_{et}$ one has

$$\Gamma(V, p_1^* F) = Mor_Y(V, Y_1 \times_X E) = Mor_Y^G(V_1 \times_X \tilde{X}, Y_1 \times_X E) =$$

$$Mor_{Y_1 \times_X \tilde{X}}^G(V \times_Y (Y_1 \times_X \tilde{X}), (Y_1 \times_X \tilde{X})^C).$$

And $\Gamma(V, p_2^* F) = Mor_{Y_2 \times_X \tilde{X}}^G(V \times_Y (Y_2 \times_X \tilde{X}), (Y_2 \times_X \tilde{X})^C)$.

Thus to construct an isomorphism $\Phi : p_1^* \xrightarrow{\cong} p_2^*$ it suffices to construct a $G$-equivariant isomorphism of $Y$-schemes

$$\phi : Y_1 \times_X \tilde{X} \xrightarrow{\cong} Y_2 \times_X \tilde{X}.$$

where $G$ acts on $Y_i \times_X \tilde{X}$ through the second factors.

Let $T$ be an $X \times X$-scheme such that $T \times_{(X \times X)} (\tilde{X} \times \tilde{X}) \cong \amalg_S \tilde{X} \times \tilde{X}$. Then $S$ is a $G \times G$-set and one says that $T$ corresponds to the $G \times G$-set $S$. Observe that in this case $T = T_S$ is an affine scheme and its coordinate ring is $k[T_S] = Map^{G \times G}(S, k[\tilde{X} \times \tilde{X}])$. A morphism $S_1 \xrightarrow{\phi} S_2$ of $G \times G$-sets induces a homomorphism of $k[X \times X]$-algebras $k[T_{S_2}] \to k[T_{S_1}]$. In particular, projections $S_1 \times S_2 \to S_i$ induce homomorphisms of $k[X \times X]$-algebras $k[T_{S_i}] \to k[T_{S_1 \times S_2}]$ and it is easy to see that the homomorphism

$$k[T_{S_1}] \otimes_{k[X \times X]} k[T_{S_2}] \to k[T_{S_1 \times S_2}]$$

is an isomorphism. Now denote by $G_1$ the set $G$ endowed with $G \times G$-action $(g_1, g_2) \cdot g = g_1 g$ and denote by $G_2$ the set $G$ endowed with $G \times G$-action $(g_1, g_2) \cdot g = g_2 g$. Let $G_{lr}$ denote the set $G$ endowed with $G \times G$-action $(g_1, g_2) \cdot g = g_1 g g_2^{-1}$. There is a commutative diagram of $G \times G$-sets

$$\begin{array}{ccc} G_{lr} \times G_1 & \xrightarrow[\cong]{\phi} & G_{lr} \times G_2 \\ {\scriptstyle pr} \downarrow & \swarrow {\scriptstyle pr} & \\ G_{lr} & & \end{array}$$

where an isomorphism $\phi$ is given by $(g_1, g_2) \mapsto (g_1, g_1^{-1} g_2)$. Observe that $\tilde{X} \times X$ corresponds to $G_1$, $X \times \tilde{X}$ corresponds to $G_2$ and $Y$ corresponds to $G_{lr}$. Thus one has a commutative diagram of $k[X \times X]$-algebras:

$$\begin{array}{ccc} k[Y] \otimes_{k[X \times X]} k[\tilde{X} \times X] & \xleftarrow[\cong]{\phi} & k[Y] \otimes_{k[X \times X]} k[X \times \tilde{X}] \\ \uparrow & \nearrow & \\ k[Y] & & \end{array}$$

and the isomorphism $\phi$ induces the desired isomorphism of $Y$-schemes. $\square$

**7.4. Corollary.** *Let $\mathcal{C}$ be a bounded complex $0 \to \mathcal{C}^{-r} \to \cdots \to \mathcal{C}^r \to 0$ with $\mathcal{C}^i \in C(\tilde{X}/X)$ then one has an isomorphism $\Phi_\mathcal{C} : p_1^* \mathcal{C} \xrightarrow{\cong} p_2^* \mathcal{C}$ of sheaf complexes on $Y_{et}$.*

*Proof.* The family $\{\Phi_{\mathcal{C}^i}\}$ of isomorphisms from 7.3 determine the desired isomorphism $\Phi_\mathcal{C}$. $\square$

**7.5. Corollary.** *Let $\mathcal{C}$ be a bounded complex with $\mathcal{C}^i \in C(\tilde{X}/X)$. Then the cohomology functor $H^p : Cp(X) \to Ab$ given by $(T, Z) \mapsto H_Z^p(T, \mathcal{C})$ is a functor with finite monodromy (see 5.3 for definition)*

*Proof.* For a couple $(T, Z) \in Cp(X)$ the isomorphism $\Phi_\mathcal{C}$ induces the isomorphism
$$H^p(\Phi_\mathcal{C}) : H_{Z_1}^p(T_1, p_1^*\mathcal{C}) \to H_{Z_2}^p(T_2, p_2^*\mathcal{C}),$$
where $Z_i = Z \times_X Y_i$ and $T_i = T \times_X Y_i$ are respective fibres. Clearly the family $\{H^p(\Phi_\mathcal{C})\}$ gives the functor isomorphism $p_1^*(H^p) \xrightarrow{\cong} p_2^*(H^p)$. $\square$

**7.6. Corollary.** *Let $\mathcal{C}$ be a bounded complex of locally constant constructible sheaves on $X_{et}$. Then the cohomology functor $H^p : Cp(X) \to Ab$ given by $(T, Z) \mapsto H_Z^p(T, \mathcal{C})$ is a functor with finite monodromy.*

*Proof.* Let $X_i \to X$ be a finite étale morphism such that $\mathcal{C}^i|_{X_i}$ is a constant sheaf on $(X_i)_{et}$. Let $\tilde{X}$ be a connected component of $X_{-r} \times_X \ldots \times_X X_r$. Then the sheaf $\mathcal{C}^i|_{\tilde{X}}$ is constant ($i = -r, \ldots, r$). Now Corollary 7.5 shows that the functor $H^p$ is a functor with finite monodromy. $\square$